\documentclass[12pt,reqno]{amsart}
\usepackage{graphicx}
\usepackage{psfrag}
\usepackage{mathrsfs}
\usepackage{color}
\usepackage{amsmath,amsthm,amsfonts,amssymb,amscd}
\font\de=cmssi12

\numberwithin{equation}{section}

\newcommand{\diff}{\operatorname{Diff}}

\newcommand{\ep} {\varepsilon} \newcommand{\eps} {\varepsilon}

\newcommand{\eqm}{\overset{\circ}{=}}

\newcommand{\ti}{\pitchfork}
\newcommand{\PW}{W}

\def \ZZ {{\mathbb Z}}

\def \RR {{\mathbb R}}

\newcommand{\PB}{\cR_{\eps,l}}

\def \cR {{\mathcal R}}
\def \cS {{\mathcal S}}

\newtheorem{theorem}{Theorem}[section]
\newtheorem{corollary}[theorem]{Corollary}
\newtheorem{lemma}[theorem]{Lemma}
\newtheorem{proposition}[theorem]{Proposition}
\newtheorem{definition}[theorem]{Definition}

\begin{document}

\thanks{ }

\author{F. Rodriguez Hertz}
\address{IMERL-Facultad de Ingenier\'\i a\\ Universidad de la
Rep\'ublica\\ CC 30 Montevideo, Uruguay.}
\email{frhertz@fing.edu.uy}\urladdr{http://www.fing.edu.uy/$\sim$frhertz}

\author{M. A. Rodriguez Hertz}
\address{IMERL-Facultad de Ingenier\'\i a\\ Universidad de la
Rep\'ublica\\ CC 30 Montevideo, Uruguay.}
\email{jana@fing.edu.uy}\urladdr{http://www.fing.edu.uy/$\sim$jana}

\author{A. Tahzibi}
\address{Departamento de Matem\'atica,
  ICMC-USP S\~{a}o Carlos, Caixa Postal 668, 13560-970 S\~{a}o
  Carlos-SP, Brazil.}
\email{tahzibi@icmc.sc.usp.br}\urladdr{http://www.icmc.sc.usp.br/$\sim$tahzibi}

\author{R. Ures}
\address{IMERL-Facultad de Ingenier\'\i a\\ Universidad de la
Rep\'ublica\\ CC 30 Montevideo, Uruguay.} \email{ures@fing.edu.uy}
\urladdr{http://www.fing.edu.uy/$\sim$ures}

\thanks{}



\renewcommand{\subjclassname}{\textup{2000} Mathematics Subject Classification}


\setcounter{tocdepth}{2}

\title[Uniqueness of SRB measures]{Uniqueness of SRB measures for transitive diffeomorphisms on surfaces}

\begin{abstract}
We give a description of ergodic components  of SRB measures in
terms of ergodic homoclinic classes associated to hyperbolic periodic points. For transitive surface
diffeomorphisms, we prove that there exists at most one SRB measure.
\end{abstract}

\maketitle

\section{Introduction}
In this paper we attempt to give a more accurate description of the ergodic components of SRB measures. These measures were introduced by Sinai, Ruelle and Bowen in the 70's (see \cite{sinai, r1, r2, b}) and are the measures most compatible with the ambient volume when the system is not conservative.

Sinai-Ruelle-Bowen's works showed the existence and some desirable  properties  of such measures for uniformly hyperbolic systems. Subsequently, SRB measures were shown to exist for many non-hyperbolic systems such as: diffeomorphisms preserving smooth measures (\cite{pesin1977}), H\'{e}non's attractors (\cite{by}), attractors with mostly contracting center direction (\cite{psinai, bv}), mostly expanding case (\cite{abv}), partially hyperbolic attractors with one-dimensional center (\cite{cy}, see also \cite{tsujii}), when $u$-Gibbs measures are unique (\cite {d1,d2}). Most of these results also include a proof of uniqueness or at least finiteness of SRB measures. In general, uniqueness results are based on the knowledge of the geometry of the unstable ``foliation". In this paper, we give a description of the ergodic components of SRB measures in terms of ergodic homoclinic classes (see next subsection) associated to periodic points. Ergodic homoclinic classes were introduced  by the authors in \cite{rhrhtu2009} (see also \cite{rhrhtu2008}) for the conservative setting. Although SRB measures have a different nature (in general, an SRB measure for $f$ is not SRB for $f^{-1}$) we obtain a similar description that combined with a subtle use of the Sard's Theorem, allows us to prove that
transitive surface diffeomorphisms have at most one SRB measure.

\subsection{Statement of results}

Roughly speaking an SRB measure is an invariant measure that has a positive Lyapunov exponent a.e. and the decomposition of the measure along unstable manifolds is equivalent to the volume. See \ref{def_srb} for a precise definition.

By Ledrappier-Young \cite{ledrappieryoung1985},  a measure satisfies the entropy formula $h_{\mu} = \int \sum_{\lambda_i > 0} \lambda_i d \mu$  if and only if it is SRB.  If all the Lyappunov exponents are non zero $\mu-$almost everywhere  then $\mu$ is called a hyperbolic measure.

 We call $\mu$ a physical measure if the basin
of $\mu, B(\mu)$ has positive Lebesgue measure, where by definition
for every continuous observable $\phi : M \rightarrow \mathbb{R}$,
$$
  \frac{1}{n} \sum_{i =0}^{n-1} \phi(f^i(x)) \rightarrow \int \phi d
  \mu
$$
for every $x \in B(\mu)$.
These measures describe the
asymptotic average behavior of a large subset of points of
the ambient space and are the basis of the understanding of
dynamics in a statistical sense.

  In general, using absolute continuity of unstable lamination for $C^{1+\alpha}$ diffeomorphisms, it turns out that any ergodic SRB measure is physical if all of its Lyapunov exponents are non-zero, see \cite{ps}. On the one hand, SRB measures are better to lead with thanks to the information given by the presence of the positive exponent. On the other hand, physical measures carry little information (see \cite{lsy_survey} for a discussion on the subject). For this reason we focus on the study of SRB measures.

  In this paper we give an accurate description for the ergodic components of  SRB measures. We define {\it ergodic homoclinic classes} for hyperbolic periodic points which are ``ergodic" version of homoclinic classes  and prove that ergodic components of hyperbolic  measures are in fact {\it ergodic homoclinic classes}.

Given a hyperbolic periodic point $p$, let us define the {\it ergodic homoclinic class} of $p, \Lambda(p),$ as the set of regular points $x \in M$ such that

\begin{figure}[h]
\begin{minipage}{7cm}
\includegraphics[width=6cm]{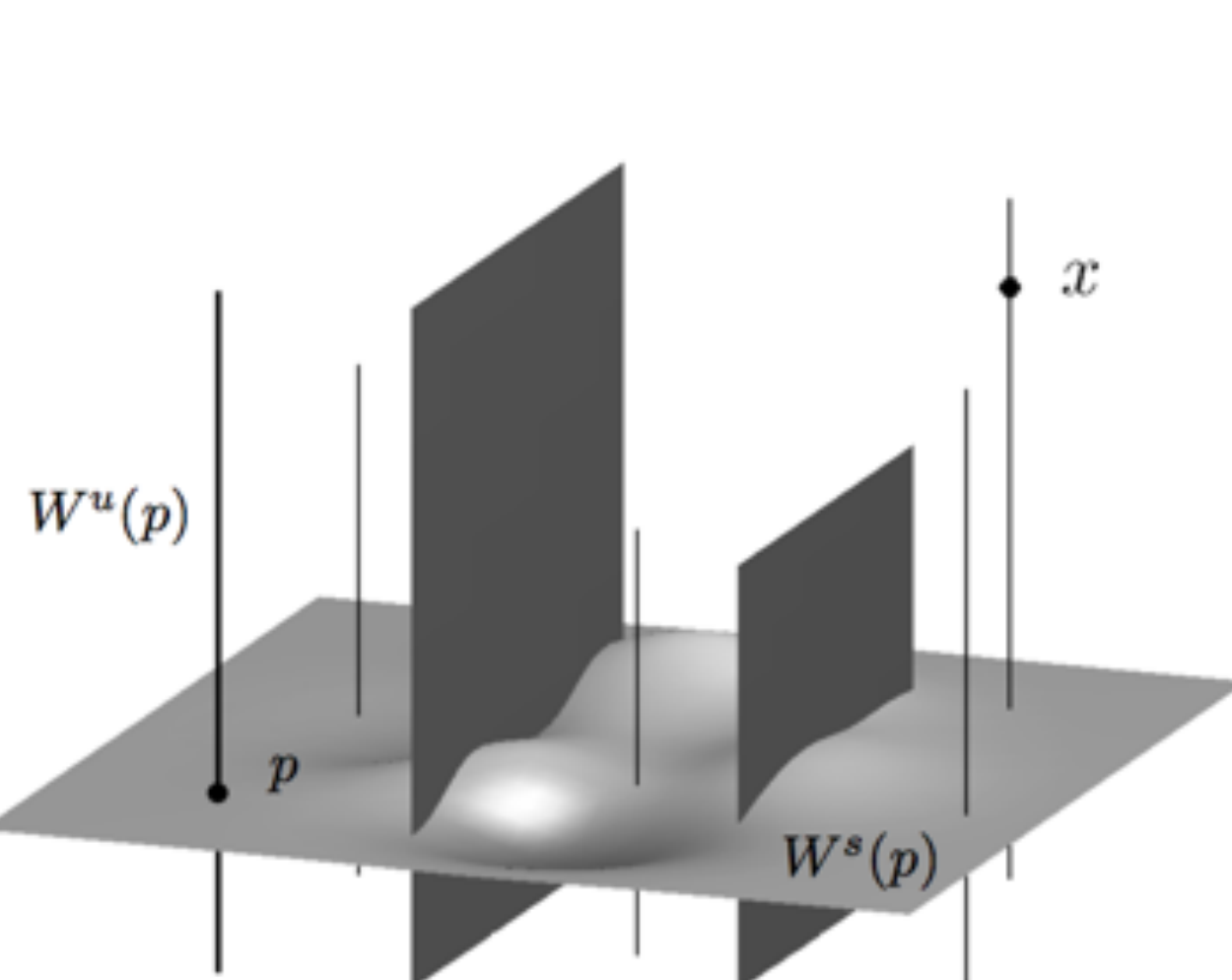}
\end{minipage}
\begin{minipage}{5cm}
\begin{eqnarray}\label{Bup}W^s(o(p))\pitchfork W^u(x)\ne\emptyset&\quad&\\
\label{Bsp}W^u(o(p))\pitchfork W^s(x)\ne\emptyset&\quad&
\end{eqnarray}
\end{minipage}
\caption{$x$ in the ergodic homoclinic class of $p$}
\end{figure}

Here $W^s(x)$ is the {\it Pesin stable manifold} of $x$, that is,
$$W^s(x)=\left\{y\in M: \lim\sup_{n\to +\infty}\frac{1}{n}\log d(f^n(x),f^n(y))<0\right\}$$
and $W^u(x)$, is {\it Pesin unstable manifold} of $x$.
For almost every point, Pesin stable and unstable manifolds are, indeed, immersed manifolds.

Note that we can write  an ergodic homoclinic class as the intersection of two invariant sets:
$$
 \Lambda(p)= \Lambda^s(p) \cap \Lambda^u(p)
$$
where $\Lambda^u(p)$ is the set of regular $x$ satisfying the relation  (1.1) and $\Lambda^s(p)$ is the set of regular points $x$ satisfying (1.2)  . It is clear that $\Lambda^s(p)$ and $\Lambda^u(p)$ are respectively $s-$saturated and $u-$saturated.

\begin{theorem}\label{criterium_srb} Let $f:M\to M$ be a $C^{1+\alpha}$ diffeomorphism over a
 compact manifold
$M$ and $\mu$ a hyperbolic SRB measure. If $\mu (\Lambda^s(p))$ and $\mu
(\Lambda^u(p))>0$, then
$$
\Lambda^u(p)\overset{\circ}{\subset} \Lambda^s(p),
$$
Moreover, the restriction of $\mu$ to $\Lambda(p)$ is ergodic and
non-uniformly hyperbolic and physical.
\end{theorem}

We can drop the hypothesis of the hyperbolicity of the measure under the hypothesis $m (\Lambda^s(p)) > 0$ where $m$ is the Lebesgue
measure. We also give an example where $\mu (\Lambda^s(p))$ and $\mu
(\Lambda^u(p))>0$ and $\Lambda^u$ is not a.e. contained in $\Lambda^s.$ Of course in such example $\mu$ is neither hyperbolic nor ergodic.

\begin{theorem} \label{criterium_srb2}
Let $f:M\to M$ be a $C^{1+\alpha}$ diffeomorphism over a compact
manifold $M$ and $\mu$ an SRB measure. If $m (\Lambda^s(p)) >0$ and $\mu
(\Lambda^u(p))>0$, then
$$
\Lambda^u(p)\overset{\circ}{\subset} \Lambda^s(p),
$$
Moreover, the restriction of $\mu$ to $\Lambda(p)$ is a hyperbolic
ergodic measure.

\end{theorem}

We mention that a similar result for Lebesgue measure has been proved
in \cite{rhrhtu2009} without the hypothesis of hyperbolicity of measure.
In fact for Lebesgue measure we did not assume the hyperbolicity the
measure and conclude that $\Lambda^s \eqm \Lambda^u.$ But, as we have mentioned before,  SRB measures have a
different nature. An SRB measure for $f$ is not SRB for
$f^{-1}.$

\par Theorem \ref{criterium_srb} has the following corollary:
\begin{corollary}\label{corollary_srb}
  Let $f:M\to M$ be a $C^{1+\alpha}$ diffeomorphism.
  If $\mu (\Lambda(p))>0$ for a hyperbolic point $p$, then $\mu |_{\Lambda(p)}$ is an ergodic component of $\mu$.
\end{corollary}

If $\mu$ is a hyperbolic invariant measure, by a result of A. Katok (\cite{katok1980}),
$(f, \mu)$ is approximated by uniformly hyperbolic (homoclinic class of hyperbolic periodic point) sets with $\mu-$measure zero. Using carefully the construction of such periodic points, we prove the following theorem which will be used in the proof of Theorem \ref{srb2}.

\begin{theorem} \label{periodic}

Let $f:M \to M$ be a $C^{1+\alpha}$ diffeomorphism over a compact
manifold $M$ and $\mu$ a hyperbolic SRB measure. Then for any
ergodic component $\nu$ of $\mu$ there exists a hyperbolic periodic
point $P$ such that $\nu (\Lambda(P)) = 1.$
\end{theorem}

For further use, we state the following simple corollary of the (inclination) $\lambda$-lemma:
\begin{proposition}\label{homoclinic}
If $p, q$ are homoclinically related then $\Lambda(p) = \Lambda(q)$.
\end{proposition}

\subsection{Uniqueness of SRB measures}

It is a challenging problem in the Ergodic Theory of
Dynamical Systems to prove the existence and uniqueness (finiteness) of SRB or physical measures.

 In this paper we  show that ergodic homoclinic classes are usefull object to distinguish between different hyperbolic ergodic SRB measures. more precisely:
 \begin{theorem} \label{uniqueness}
 Let $\mu, \nu$ be ergodic SRB measures such that  $\mu(\Lambda(p)) = \nu(\Lambda(p)) =1$ for some hyperbolic periodic point $p$ then $\mu= \nu.$
 \end{theorem}

 In the two dimensional case, using the above theorem we  prove that topological transitivity is enough to guarantee that there exists at most one SRB measure for surface $C^2$- diffeomorphisms.

 It is easy to see that for a surface diffeomorphism any
SRB measure is hyperbolic. Indeed take an ergodic SRB measure $\mu$
with two Lyapunov exponents $\lambda^+
> \lambda^{-}.$ By Pesin entropy formula (see Young-Ledrappier \cite{ledrappieryoung1985}) $h(\mu)
= \lambda^{+}$. Since  $h(\mu, f) = h(\mu, f^{-1})$ by Ruelle inequality
$\lambda^{+} \leq - \lambda^{-}$.

\begin{theorem} \label{srb2}
Let $f:M\to M$ be a $C^{1+\alpha}$ diffeomorphism over a compact
surface $M$. If $f$ is topologically transitive then there exist at
most one SRB measure.
\end{theorem}
An example due to I. Kan (see \cite{kan}) shows that the above theorem
can not be true in higher dimensional manifolds. I. Kan constructed a
transitive diffeomorphism of $\mathbb{T}^2\times [0,1]$ with two SRB measures with
intermingled basins. Gluing along the boundary torus two of these examples
and composing with a diffeomorphism that interchanges the two components a transitive diffeomorphism of $\mathbb{T}^3$ with two SRB-measures is obtained.

 We recall that,  as a consequence of the absolute continuity of the unstable ``foliation", Lebesgue measure is SRB if there is a positive exponent a.e. (see \cite{pesin1977}). As a corollary of the above theorem we obtain the following result in the conservative setting. We thank F. Ledrappier for observing this point in the Workshop on Partial Hyperbolicity, in Beijing.

\begin{theorem}
Let $f: M \to M$ be a $C^{1+\alpha}$ volume preserving diffeomorphism of a compact surface $M$ with non zero Lyapunov exponents.  If $f$ is topologically transitive then it is ergodic.
\end{theorem}

Observe that by H. Furstenberg's example (\cite{furst}), without the hypothesis of non-vanishing Lyapunov exponents, even minimality is not enough to guarantee ergodicity. Furstenberg constructs a minimal non-ergodic $C^\infty$ diffeomorphism of the two torus.

\section{Preliminaries}\label{section.preliminaries}

\subsection{Non-uniform hyperbolicity}\label{subsection.pesin}
Let us review some results about Pesin theory that shall be used in this paper. A good summary of these facts may be found, for instance, in
\cite{ps} and \cite{ledrappieryoung1985}. For further references, see
 Katok's paper \cite{katok1980} and  the book by Barreira and Pesin \cite{barreira-pesin}.\par
Let $f:M\to M$ be a $C^1$ diffeomorphism  of a compact Riemannian manifold of dimension $n$.
Given a vector $v\in T_xM$, let the {\de Lyapunov exponent of $v$} be the exponential growth rate of $Df$ along $v$, that is
\begin{equation}\label{lyap exp}
\lambda(x,v)=\lim_{|n|\to\infty}\frac1n\log|Df^n(x)v|
\end{equation}
in case this amount is well defined. And let $E_\lambda(x)$ be the subspace of $T_xM$ consisting of all $v$ such that the Lyapunov exponent of $v$ is $\lambda$. Then we have the following:
\begin{theorem}[Osedelec]
For any $C^1$ diffeomorphism $f:M\to M$ there is an $f$-invariant Borel set $\cR$ of total probability (in the sense that $\mu(\cR)=1$ for all invariant probability measures $\mu$), and for each $\eps>0$ a Borel function $C_\eps: \cR\to(1,\infty)$
such that for all $x\in \cR$, $v\in T_xM$ and $n\in\ZZ$
\begin{enumerate}
\item $T_xM=\bigoplus_\lambda E_\lambda(x)$ ({\de Oseledec's splitting})
\item For all $v\in E_\lambda(x)$
$$C_\eps(x)^{-1}exp[(\lambda-\eps)n]|v|\leq|Df^n(x)v|\leq C_\eps(x)exp[(\lambda+\eps)n]|v|$$
\item $\angle\left(E_\lambda(x),E_{\lambda'}(x)\right)\geq C_\eps(x)^{-1}$ if $\lambda\ne\lambda'$
\item $C_\eps(f(x))\leq \exp(\eps) C_\eps(x)$
\end{enumerate}
\end{theorem}
The set $\cR$ is called the set of {\de regular points}.
We also have that $Df(x)E_\lambda(x)=E_\lambda(f(x))$. If
an $f$-invariant measure $\mu$ is ergodic then the Lyapunov exponents and $\dim E_\lambda(x)$ are
constant $\mu$-a.e.\par
For fixed $\ep > 0$ and given $l > 0,$ we define the {\de Pesin blocks}:
 $$
  \PB = \left\{ x \in \mathcal{R} : C_\eps(x) \leq l\right\}.
 $$
Note that Pesin blocks are not necessarily invariant. However $ f(\PB) \subset \cR_{ \eps,\exp(\eps)l}$.
Also, for each $\eps>0$, we have
 \begin{equation}\label{pesin region}
\cR = \bigcup_{ l= 1}^{\infty} \PB
\end{equation}
We loose no generality in assuming that $\PB$ are compact.
For all $x\in \cR$ we have
$$
T_xM=\bigoplus_{\lambda<0} E_\lambda(x)\oplus E^0(x)\bigoplus_{\lambda>0} E_\lambda(x)
$$
where $E^0(x)$ is the subspace generated by the vectors having zero Lyapunov exponents. Let $\mu$ be an invariant measure. When $E^0(x)=\{0\}$ for $\mu$-a.e. $x$ in a set $N$, then we say that $f$ is {\de non-uniformly hyperbolic} on $N$ and that $\mu$ is a {\de hyperbolic measure} on $N$.\par
Now, let us assume that $f\in C^{1+\alpha}$ for some $\alpha>0$. Given a regular point $x$, we define its
{\de stable Pesin manifold} by
\begin{equation}\label{stable Pesin manifold}
\PW^s(x)=\left\{y:\limsup_{n\to +\infty}
\frac{1}{n}\log  d(f^n(x),f^n(y))<0\right\}
\end{equation}
The {\de unstable Pesin manifold} of $x$, $\PW^u(x)$ is the stable Pesin manifold of $x$ with respect to $f^{-1}$.
Stable and unstable Pesin manifolds of points in $\cR$ are immersed manifolds \cite{pesin1977}. We stress that $C^{1+\alpha}$ regularity is crucial for this to happen. In this way we obtain a partition
 $x\mapsto\PW^s(x)$ , which we call {\de stable partition}. Unstable partition is defined analogously. Stable and unstable partitions are invariant.\par
On the Pesin blocks we have a continuous variation: Let us call $W^s_{loc}(x)$ the connected component of $W^s(x)\cap B_r(x)$ containing $x$, where $B_r(x)$ denotes the Riemannian ball of center $x$ and radius $r>0$, which is sufficiently small but fixed. Then
\begin{theorem}
[Stable Pesin Manifold Theorem \cite{pesin1977}]
Let $f:M\to M$ be a $C^{1+\alpha}$ diffeomorphism preserving a smooth measure $m$. Then, for each $l>1$ and small $\eps>0$, if $x\in \cR_{\eps,l}$:
\begin{enumerate}
\item $W^s_{loc}(x)$ is a disk such that $T_xW^s_{loc}(x)=\bigoplus_{\lambda<0} E_\lambda(x)$
\item $x\mapsto\PW^s_{loc}(x)$ is continuous over $\PB$ in the $C^1$ topology
\end{enumerate}
\end{theorem}
In particular, the dimension of the disk $W^s_{loc}(x)$ equals the number of negative Lyapunov exponents of $x$. An analogous statement holds for the unstable Pesin manifold.
\subsection{Absolute continuity}
An important notion behind the criterion we are going to prove is absolute continuity. Let us state the definitions we will be using. The point of view we follow is similar to that in \cite{ledrappieryoung1985}.\par
Let $\xi$ be a partition of the manifold $M$. We shall call $\xi$ a {\de measurable partition} if the quotient space $M/\xi$ is separated by a countable number of measurable sets. For instance, the partition of the 2-torus by lines of irrational slope is not measurable, while the partition of $[0,1]$ by singletons is measurable. The quotient space $M/\xi$ of a Lebesgue space $M$ by a measurable partition $\xi$ is again a Lebesgue space \cite{rohlin}.\par
Associated to each measurable partition $\xi$ of a Lebesgue space $(M,{\mathcal B},m)$ there is a canonical system of {\de conditional measures} $m_x^\xi$, which are measures on $\xi(x)$, the element of $\xi$ containing $x$, and with the property that for each $A\in{\mathcal B}$ the set $A\cap\xi(x)$ is measurable in $\xi(x)$ for almost all $\xi(x)$ in $M/\xi$, and the function $x\mapsto m_x^\xi(A\cap\xi(x))$ is measurable, with:
\begin{equation}\label{ecuacion.measurable.partition}m(A)=\int_{M/\xi}m_x^\xi(A\cap\xi(x))dm_T\end{equation}
where $m_T$ is the quotient measure on $M/\xi$. For each measurable partition this canonical system of conditional measures is unique (mod $0$), i.e. any other system is the same for almost all $\xi(x)\in M/\xi$. Conversely, if there is a canonical system for a partition, then the partition is measurable. In our case, we will be interested in stable and unstable partitions, note that in general these partitions are not measurable.\par
A measurable partition $\xi$ is {\de subordinate} to the unstable partition $W^u$ if for $m$-a.e. we have $\xi(x)\subset W^u(x)$, and $\xi(x)$ contains a neighborhood of $x$ which is open in the topology of $W^u(x)$.
\begin{definition}\label{fubini.like}
  $m$ has {\de absolutely continuous conditional measures on unstable manifolds} if for every measurable partition $\xi$ subordinate to $W^u$, $m^\xi_x<<\lambda^u_x$ for $m$-a.e. $x$, where $\lambda^u_x$ is the Riemannian measure on $W^u(x)$ given by the Riemannian structure of $W^u(x)$ inherited from $M$.
\end{definition}

We are now able to give a definition of SRB-measure.

\begin{definition}\label{def_srb}
An $f$-invariant probability measure $\mu$ is called a Sinai-Bowen Ruelle (SRB) measure if it has a positive Lyapunov exponent a.e. and absolutely continuous conditional measures on unstable manifolds. After Ledrappier-Young (\cite{ledrappieryoung1985}) this is equivalent to having a positive Lyapunov exponent a.e. and satisfying the Pesin formula, $h_\mu(f)=\int\sum_{\lambda(x)>0}\lambda(x)\, d\mu$.
\end{definition}

Now, take a point $x_0\in \cR$, the set of regular points. Assume that $x_0$ has at least a negative Lyapunov exponent. Take two small discs $T$ and $T'$ near $x_0$ which are transverse to $W^s(x_0)$.
Then we can define the {\de holonomy map} with respect to these transversals as a map $h$ defined on a subset of $T$ such that $h(x)=W^s_{loc}(x)\cap T'$.
The domain of $h$ consists of the points $x\in T\cap \cR$ whose stable manifold have the same dimension as $W^s(x_0)$, and which transversely intersect $T$ and $T'$. $h$ is a bijection.
\begin{definition}\label{absolute.continuity}
We say that the stable partition is {\de absolutely continuous} if all holonomy maps are measurable and take Lebesgue zero sets of $T$ into Lebesgue zero sets of $T'$.
\end{definition}
Absolute continuity of the unstable partition is defined analogously.
\begin{theorem}[\cite{pesin1977}] \label{teorema.cont.abs.Pesin}
Let $f$ be a $C^{1+\alpha}$ diffeomorphism. Then, its stable and unstable partitions are absolutely continuous.
\end{theorem}
Note that the holonomy maps of the stable foliation are continuous and have continuous Jacobians when restricted to the Pesin blocks $\PB$.


\section{Ergodic components and ergodic homoclinic classes}
In this section we prove all theorems except Theorem \ref{srb2}.
For the sake of simplicity firstly, we shall prove Corollary
\ref{corollary_srb}. Let us introduce some lemmas before
entering into its proof.be an SRB-measure. Let $\mu$  For any given function $\varphi\in
L^1_{\mu}(M,\RR)$, let
\begin{equation}  \label{birkhoff_srb}
\tilde\varphi^\pm(x)=\lim_{n\rightarrow\pm
\infty}\frac1n\sum^{n-1}_{i=0}\varphi(f^n(x))
\end{equation}
By Birkhoff Ergodic Theorem, the limit (\ref{birkhoff_srb}) exists
and $\varphi^+(x)=\varphi^-(x)$ for $\mu-$a.e. $x\in M$. Note that
$\varphi^\pm(x)$ is $f$-invariant. Moreover, we have the following:
\begin{lemma}  \label{full.continuous_srb}
For all $\varphi\in C^0(M)$ there exists an invariant set $\cS$ with
$\mu(\cS)=1$ such that if $x\in \cS$ we have
$\tilde\varphi_+(w)=\tilde\varphi_+(x)$ for all $w\in W^s(x)$ and
$m^u_x$-a.e. $w\in W^u(x)$.
\end{lemma}
\begin{proof}  [Proof of the Lemma] The proof is completely
analogous to the one for smooth measures since SRB-measures, by
definition, have absolutely continuous conditional measures along
unstable manifolds.
\end{proof}

\begin{proof}[Proof of Corollary \ref{corollary_srb}]
Let $\varphi:M\to \RR$ be
a continuous function, let $\cS$ be the set obtained in Lemma \ref{full.continuous_srb} and
$\cR$ the set of $\mu$-regular points. We shall see that $\tilde\varphi^+$ is
constant on $\Lambda(p)\cap \cS\cap\cR$. This will prove that $\mu$ is
ergodic when restricted to $\Lambda(p)$.\par Let $x,y\in
\Lambda(p)\cap\cS\cap\PB=\Lambda$ for some $\eps>0$ and $l>1$. Without
loss of generality we may assume that $x$ and $y$ are in the support
of the restriction of $\mu$ to $\Lambda$ and that they return
infinitely many times to $\Lambda$. Hence there exists $n>0$ such
that $f^n(y)\in \Lambda$ and $d(f^n(y),W^u(p))<\delta/2$ where
$\delta>0$ is as in the definition of transverse absolute
continuity. Hence $W^s_{loc}(f^n(y))\ti W^u(p) \neq \emptyset$. We can suppose for
simplicity that $n=0$, and that $p$ is a hyperbolic  fixed
point.\par As a consequence of the Inclination Lemma, there exists
$k>0$ such that $f^k(x)\in\Lambda$ and $W^u(f^k(x))\ti
W^s_{loc}(y) \neq \emptyset$. As in the case of smooth measures, due to Lemma
\ref{full.continuous_srb} above, there is a $m^u_y$-positive measure
set of $w\in W^u_{loc}(y)$ such that
$\tilde\varphi^+(w)=\tilde\varphi^+(y)$ and $W^s_{loc}(w)\ti
W^u(f^k(x)) \neq \emptyset$. Since $\tilde\varphi^+$ is constant on stable leaves
and there is transverse absolute continuity, we get a
$m^u_{f^k(x)}$-positive measure set of $w'\in W^u_{loc}(f^k(x))$
such that $\tilde\varphi^+(w')=\tilde\varphi^+(y)$. See Figure
\ref{figure.proofbp}. But $f^k(x)\in\cS$, so due to Lemma
\ref{full.continuous_srb} again
$\tilde\varphi^+(y)=\tilde\varphi^+(f^k(x))=\tilde\varphi^+(x)$
concluding  the proof.
\par

\end{proof}
In order to prove Theorem \ref{criterium_srb}, we shall need a
refinement of Lemma \ref{full.continuous_srb}:
\begin{lemma} \label{fullmeasure_srb} Given $\varphi\in L^1(\mu)$ there exists an invariant set $\cS_{\varphi}\subset M$,
$\mu(\cS_{\varphi})=1$ such that if $x\in\cS_{\varphi}$ then $m^u_x$-a.e. $y\in W^u(x)$ satisfy $\varphi^+(y)=\varphi^+(x)$.
\end{lemma}
\begin{proof}
Given $\varphi\in L^1(\mu)$ take a sequence of continuous functions
$\varphi_n$  converging to $\varphi$ in $L^1(\mu)$. Since
$\tilde\varphi^+_{n}$ converges in $L^1(\mu)$ to $\tilde\varphi^+$
there exists a subsequence $\tilde\varphi_{n_k+}$ converging a.e. to
$\tilde\varphi_+$. The intersection of this set of almost every
where convergence with the set $\cS$ obtained in Lemma
\ref{full.continuous_srb} gives the desired set $\cS_\varphi$.
 \end{proof}

 Let us give now the proof of Theorems \ref{criterium_srb} and \ref{criterium_srb2}.

\begin{proof}[Proof of Theorem \ref{criterium_srb}]To simplify ideas, let us suppose that $p$ is a hyperbolic fixed point.
Let $\cS$ be the set obtained in Lemma \ref{fullmeasure_srb} for the
characteristic function ${\mathbf 1}_{\Lambda^s(p)}$.
Take $x\in \Lambda^u(p)\cap \cS$, we will prove that $x\in \Lambda^s(p)$.
This proves the first claim of the Theorem. Let $y\in \Lambda^s(p)$, and $\eps>0$, $l>1$
be such that $x,y\in\PB$. We loose no generality in assuming that
$y$ is in the support of $\mu$ restricted to $\Lambda^s(p)\cap \PB$. We
can also assume that $x$ and $y$ return infinitely many times to
$\PB$. Proceeding as in the proof of Corollary \ref{corollary_srb},
we may also assume that $d(y,W^u(p))<\delta/2$.\par Note that
${\mathbf 1}_{\Lambda^s(p)}$ is an $f$-invariant function. This implies
that if $x\notin \Lambda^s(p)$ then $m^u_x$-a.e. $y\in W^u(x)$ will
satisfy $y\notin \Lambda^s(p)$, due to Lemma \ref{fullmeasure_srb} above.
The idea of the proof is to find a $m^u_x$-positive measure set of
points $z\in W^u(x)$ such that $z\in \Lambda^s(p)$. This will prove that
$x\in \Lambda^s(p)$. \par As a consequence of the Inclination Lemma, and
since $x$ returns infinitely many times to $\PB$, there exists $k>0$
such that $W^u(f^k(x))\ti W^s_{loc}(y) \neq \emptyset$. Note that this intersection a priori can have positive dimension.\par

\begin{figure}[h]\label{figure.proofbp}
\includegraphics[width=6cm]{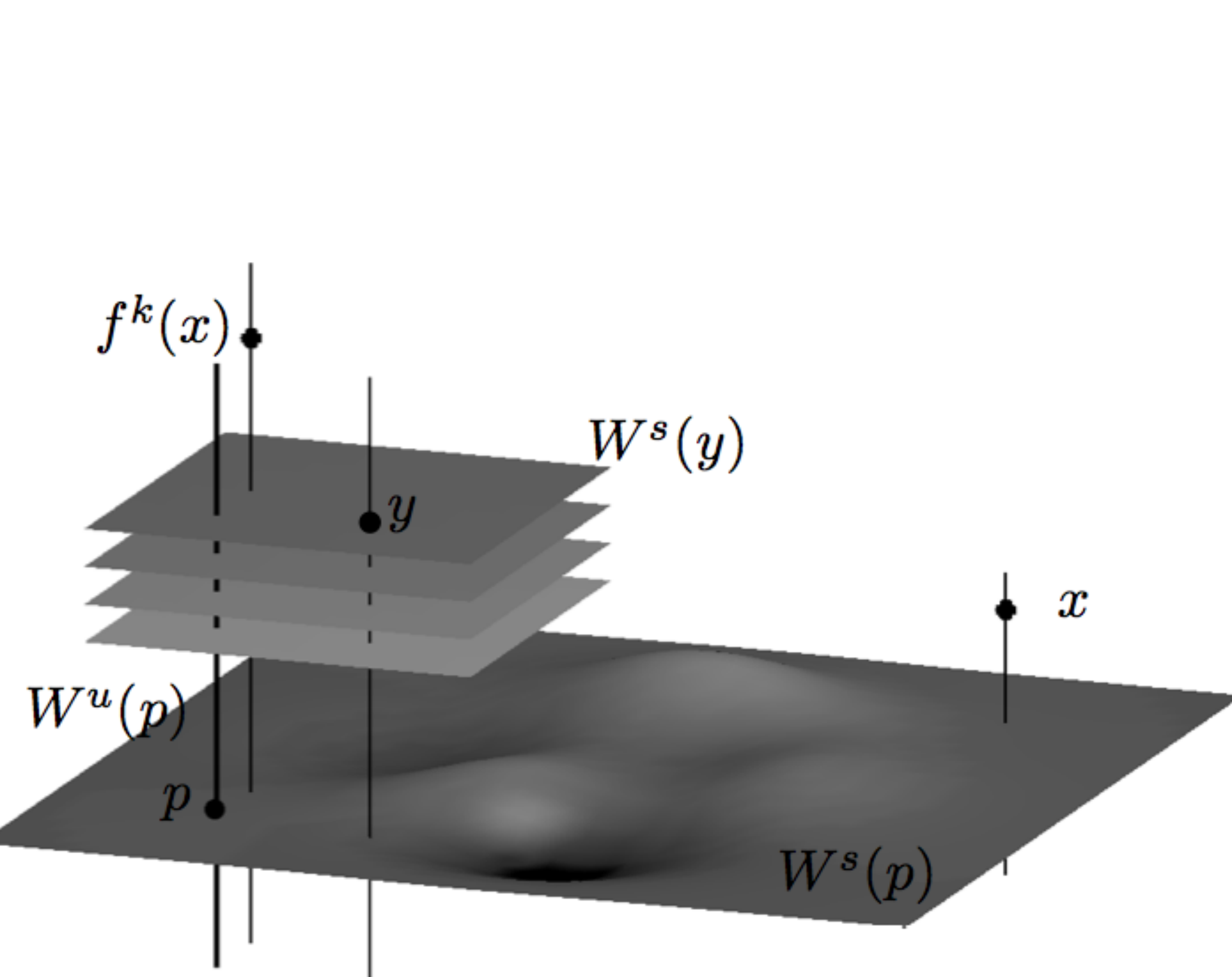}
\caption{Proof of Corollary \ref{corollary_srb}}
\end{figure}

Now, since $y$ is in the support of $\mu$ restricted to $\PB\cap
\Lambda^s(p)$, we have $\mu(\PB\cap \Lambda^s(p)\cap B_\delta(y))>0$.
 Since by hypothesis $\mu$ is a hyperbolic measure with absolutely continuous conditional measures along unstable manifolds (SRB-measure) there exists $z\in B_\delta(y)$ such that $\dim(W^u(y)) = \dim(W^u(z))$ and $m^u_z(W^u_{loc}(z)\cap \PB\cap \Lambda^s(p)\cap B_\delta(y)) > 0$. Take a smooth foliation $\mathcal{L}$ of a neighborhood of a point of $W^s_{loc}(y)\cap W^u(f^k(x))$ inside $W^u(f^k(x))$ of dimension equal to $\dim(W^u(y)) = n - \dim(W^s(y))$. This can be done in such a way that every $L\in \mathcal{L}$ is transversal to $W^s_{loc}(y)$. In fact
 \begin{align*}
 \dim(W^s(y) \cap W^u(f^k(x)) &= \dim(W^u(f^k(x))) + \dim(W^s(y)) - n \\
 &= \dim(W^u(f^k(x))) - \dim(W^u(y)).
 \end{align*}

Take an open submanifold $T$ of $W^u(f^k(x)) \ti W^s_{loc}(y).$ Transverse  absolute continuity of the stable foliation implies that
$m^L_\omega(\Lambda^s(p)\cap L)>0$ $\forall \omega \in T$ where $m^{L}$ is the induced Lebesgue measure of the the leaf $L$ of the smooth foliation $\mathcal{L}$. Now Fubini theorem for the smooth foliation $\mathcal{L}$ implies that
\begin{equation*}m^u_{f^k(x)}(\Lambda^s(p)\cap W^u(f^k(x))) \geq \int_{T} m^{L}_{\omega} (L_{\omega} \cap \Lambda^s(p)) dm_T(\omega) > 0
\end{equation*}
So we have proved that a $m_{f^k(x)}^u$-positive measure subset of $W^u(f^k(x))$ belongs to $\Lambda^s(p).$ As $x \in \mathcal{S}$ and $\mathcal{S}$ is invariant this implies that $f^k(x) \in \Lambda^s(p)$ and so $x \in \Lambda^s(p).$ This finishes the proof of $\Lambda^u(p)\overset{\circ}{\subset} \Lambda^s(p)$.

By a similar argument as above we are able to prove the ergodicity of $\mu|_{\Lambda(p)}.$ Indeed, let $\varphi \in C(M)$ and $x, y \in \Lambda(p)$ be as above with $\cS$ be the full measure obtained in lemma \ref{full.continuous_srb} for $\varphi.$ Following the above arguments mutatis mutandis, we prove that $\varphi_{+} (x) = \varphi_{+}(y).$

By definition the restriction of $f$ on $\Lambda(p)$ is non uniformly hyperbolic with the same index of $p.$

\end{proof}

%

\begin{proof}[Proof of Theorem \ref{criterium_srb2}]
  Let $x, y$ as in the proof of Theorem \ref{criterium_srb} with the difference that now $y$ is a density point of $\PB\cap
\Lambda^s(p)$ for Lebesgue measure.  Now take an smooth foliation $\mathcal{F}$ of dimension equal to $\dim(W^u(y))$ inside $B_{\delta}(y).$ By Fubini theorem there exist a leaf of this foliation $\mathcal{F}(z), z \in B_{\delta}(y)$ such that   $m^{\mathcal{F}}_z( \mathcal{F}(z)\cap \PB\cap \Lambda^s(p)\cap B_\delta(y)) > 0$
 where here $m_z^{\mathcal{F}}$ denotes the Lebesgue measure of $\mathcal{F}(z).$ From this point, just changing the role of $W^u(z)$ by $\mathcal{F}(z)$ in the proof of theorem  \ref{criterium_srb}  the  arguments are the same.
 \end{proof}

\begin{proof}[Proof of Theorem \ref{periodic}]
Let $\mu$ be a hyperbolic SRB measure. Firstly we prove the following well-known fact.

\begin{lemma} Almost all ergodic components of $\mu$ are hyperbolic and SRB.
\end{lemma}
\begin{proof}
By ergodic decomposition there exists a probability measure $\hat{\mu}$ on $\mathcal{M}(\mathcal{M}(f))$ with support on ergodic measures such that $$h_{\mu} = \int_{\mathcal{M}(f)} h_{\nu} d\hat{\mu}(\nu) ,$$
$$
\int \sum \lambda^{+}_i d \mu = \int \sum \lambda^+_i (\nu) d {\hat{\mu}}(\nu).$$

By Ruelle inequality we have that for all $\nu, h(\nu) \leq \int \sum \lambda^+_i (\nu)$  and putting these together it is clear that  $\hat{\mu}-$almost every $\nu$ will satisfy the entropy formula and so it is an SRB measure.
\end{proof}

From now on suppose that $\mu$ is itself an ergodic hyperbolic SRB measure. By a Katok result
$(f, \mu)$ is approximated by uniformly hyperbolic sets with $\mu-$measure zero. More precisely, take a Pesin block $\Gamma$ of large measure such that the size of stable and unstable manifolds are uniformly bounded from below by a constant larger than zero.
 Let $x \in supp(\mu | \Gamma)$ and $B$ be a small ball around $x$ such that $\mu(B \cap \Gamma) > 0.$   By Katok closing lemma we can find a periodic point $p$ near enough to $x$ whose stable and unstable manifolds respectively are $C^1-$close to the stable and unstable manifolds of $p$ and consequently $x \in \Lambda(p).$
 As the stable and unstable lamination vary continuously on $\Gamma$ we obtain that $y \in \Lambda(p)$ for any $y \in B \cap \Gamma.$ This yields that $\mu(\Lambda(p)) > 0.$ The ergodicity of $\mu$ implies that  $\mu(\Lambda(p)) =1.$
\end{proof}

\section{Example}
Here we give some examples of systems with SRB measures which shed light on the difference between our results in smooth and SRB measures case.

We construct a diffeomorphism $f : M \rightarrow M$ with an SRB measure $\mu$ such that there exists a hyperbolic periodic point $p$ such that $\mu(\Lambda^s(p)) = \mu(\Lambda^u(p)) = 1/2$ but $\mu(\Lambda^s \cap \Lambda^u)=0 .$  Observe that $\mu$ can not be the Lebesgue measure by our previous work in \cite{rhrhtu2009}. In the smooth measure case $\mu(\Lambda^s), \mu(\Lambda^u) > 0$ implies that $\Lambda^s \eqm \Lambda^u \eqm \Lambda$ and $f|\Lambda$ is ergodic.

We will split the construction of $f$ into four steps.

\subsection{First Step} We begin with $f_0 : \mathbf{T}^2 \to \mathbf{T}^2$  be a $C^{1 +\alpha},  0 < \alpha < 1$  almost -Anosov diffeomorphism with  an SRB measure $\mu_0$. Moreover $f_0$ has a fixed point $R$ such that $Df_0(R)$ has two eigenvalues $\lambda_1= 1, \lambda_2 < 1$. Such $f_0$ can be obtained satisfying the following properties:
\begin{itemize}
\item $\mu_0$-almost every $x \in \mathbf{T}^2$ has one positive (and consequently one negative) Lyapunov exponent,
\item $W^s(R) \ti W^u(P)$ for any periodic point $P \neq R,$
\item $W^u(Q) \ti W^s(P)$ for any two periodic points $P, Q \in \mathbf{T}^2.$
\end{itemize}
We emphasize that such example can not be $C^2.$ See the  work of  Hatomoto \cite{hatomoto} (see also \cite{huyoung})

\subsection{Second Step} Now we consider a family of skew products over $f_0$ as follows.  Recall that $\mu_0$ is an SRB measure of $f_0$ and $f_0$ has a fixed point $R$ with a neutral direction. We assume also that $f_0$ has  two more fixed points $P, Q$ which are hyperbolic with one dimensional unstable manifold. For $x \in \mathbf{T}^2$ and $t\geq 1$ let $g^{t}_x : \mathbb{S}^1 \rightarrow \mathbb{S}^1$  satisfy the following properties:
\begin{itemize}
\item For all $x \in \mathbf{T}^2$, $g_x^t : \mathbb{S}^1 \rightarrow \mathbb{S}^1, g_x^t (0) = 0 , | Dg_x^t (0) | \leq t,$
\item For some small $\epsilon > 0,  | Dg_x^t (0) | = t ,  \forall x \in B_{\epsilon}(Q),$
\item  $\frac{1}{3} \leq  |Dg_x^t (0)| \leq \frac{1}{2}, \forall x \notin B_{2\epsilon}(Q),$
\item $t \rightarrow \int_{\mathbf{T}^2} \log |D g^t_x (0) | d \mu_0$ is continuous.
\end{itemize}
Now we consider the following skew product over $f_0:$ $$ F^{t} (x, \theta) = (f_0(x), g^{t}_{x} (\theta) )  $$

\begin{lemma}
There exists $t_0$ such that for $\mu_0 \times \delta_0$ almost every $(x, 0) \in \mathbf{T}^2 \times \mathbb{S}^1,$ the Lyapunov exponent  of  $F^{t_0}$ in the tangent direction to $\mathbb{S}^1$ vanishes.
\end{lemma}
\begin{proof}
 Since $\mu_0$ is ergodic, let $\alpha:= \mu_0(B_{\epsilon} (Q)$.  By Birkhoff's Theorem for a  $\mu_0-$typical point $x:$
$$
 \alpha \log(t) - (1 -  \alpha) \log(3)  \leq \int_{\mathbf{T}^2} \log |D g^t_x (0) | d \mu_0 =  \lim_{n \rightarrow \infty} \frac{1}{n} \log | \Pi_{i=0}^{n-1} Dg_{f_0^{i}(x)} (0) | \leq  \log(t)$$
 Using the above estimates and the continuity of $\int_{\mathbf{T}^2} \log |D g^t_x (0) | d \mu_0$ with respect to $t$,  we conclude that there exist $t_0$ such that $ \int_{\mathbf{T}^2} \log |D g^t_x (0) | d \mu_0 = 0$ which means that  $\mu_0 \times \delta_0$-almost every point $(x , \theta) \in \mathbf{T}^2 \times  \mathbb{S}^1$ has zero Lyapunov exponent.

\end{proof}

\subsection{Third Step} Let  $F := F^{t_0}$ with $t_0$ as above. Let $A : \mathbb{T}^2 \rightarrow \mathbb{T}^2$ be a linear Anosov  diffeomorphism and define
$$
f \in \diff^{1 + \alpha}( \mathbf{T}^2 \times \mathbb{S}^1 \times \mathbb{T}^2) , \, f(x, \theta, y) = (F(x, \theta), A(y)).$$

The Lebesgue measure of $\mathbb{T}^2$ is an SRB measure for $A$ and we denote it by $m.$ So, the probability measure $ \mu := \frac{\mu_0 + \delta_{R} }{2}\times \delta_0 \times m$ is invariant by $f$. We will show that $\mu$ is SRB and satisfies $ \mu(\Lambda^s(P)) = \mu(\Lambda^u(P)) = \frac{1}{2}$ and $\Lambda^s(P) \cap \Lambda^u(P) = \emptyset.$

The SRB property is straightforward from the definition and the fact that the Lyapunov exponent of $f$ along the tangent direction to $\mathbb{S}^1$ vanishes and $\mu_0$ and $m$ are SRB.

By our construction $P$ is a hyperbolic fixed point with unstable dimension two. $Q$ is also hyperbolic but with unstable dimension three and $R$ has a neutral direction and one dimensional unstable manifold.  As $ W^s(R, f_0) \ti W^u(P, f_0)$  and $A$ is Anosov, we conclude that $ \{R\} \times \{0\} \times \mathbb{T}^2 \subset \Lambda^s$ and consequently $\mu(\Lambda^s (P) \geq \mu( \{R\} \times \{0\} \times \mathbb{T}^2) = \frac{1}{2}.$ In fact, as for $\mu_0 \times \delta_0$-almost every $(x, \theta)$ the Pesin stable manifold is two dimensional, no such point belongs to $\Lambda^s(P)$ and we hav eproved that $ \mu(\Lambda^s) = \frac{1}{2}.$

Let us now investigate $\Lambda^u(P).$ By construction $\mu_0 \times \delta_{R} \times m$-almost avery point has two dimensional unstable anifold which is transverse to $W^s(P).$ It is clear that $ \{R\} \times \{0\} \times \mathbb{T}^2 \notin \Lambda^u(P)$. We proved that $\mu(\Lambda^u(P)) = \frac{1}{2}$ and $\Lambda^s(P) \cap \Lambda^u(P) = \emptyset.$

\section{SRB measure for surface diffeomorphisms}
In this section we prove Theorem \ref{srb2}. For this aim we will first prove Theorem \ref{uniqueness}. This will be done in the next subsection.

\subsection{SRB measures supported on the same ergodic homoclinic class}
In this subsection we will show that ergodic homoclinic classes support at most one SRB measure. This is the statement of Theorem \ref{uniqueness}.

\begin{proof}[Proof of Theorem \ref{uniqueness}]

 Let $B(\mu)$ and $B(\nu)$ be, respectively, the basins of $\mu$ and $\nu$. By ergodicity we have that $\mu(B(\mu)) = \nu(B(\nu)) = 1.$ In fact  by Birkhoff's Ergodic Theorem   there  exists $B_{\mu} \subset B(\mu)$ and $B_{\nu} \subset B(\nu)$ such that $ \mu(B_{\mu}) = \nu(B_{\nu}) =1$ where
$$
 B_{\mu} = \{ x \colon \lim_{n \rightarrow  \infty}\frac{1}{n} \sum_{i =0}^{n-1} \varphi(f^{\pm i}(x)) \rightarrow \int \varphi d
  \mu  \quad \forall \varphi \in C(M) \}
$$
$$
B_{\nu} = \{ x \colon \lim_{n \rightarrow  \infty}\frac{1}{n} \sum_{i =0}^{n-1} \varphi(f^{\pm i}(x)) \rightarrow \int \varphi d
  \nu  \quad \forall \varphi \in C(M) \}
  $$

 Since $\mu(\Lambda)=\nu(\Lambda) =1$ it comes out that $\mu((B_{\mu}) \cap \Lambda) = \nu((B_{\nu}) \cap \Lambda) =1.$
From now on the technique will be similar to one of the proof of Corollary \ref{corollary_srb}.

By absolute continuity of unstable lamination we can take $x$ such that $m_x^u (B_{\mu} \cap \Lambda) = 1.$  Now let $y$ be a point of recurrence of $B_{\nu} \cap \Lambda_{\epsilon}$ where $\Lambda_{\epsilon}$ is a Pesin block for $\nu | \Lambda$. Again by absolute continuity of unstable lamination we can choose $y$ in such a way that $m^u_{y} (\Lambda_{\epsilon} \cap B_{\nu} ) > 0.$  Since $y$ returns infinitely many times to the Pesin's block we additionally can assume that $y$ is close enough to $W^u(p)$. Using  the $\lambda$-lemma we have that  $W^u(f^k(x))$ is also $C^1-$close to $W^u(p)$. This implies that $W_{\epsilon}^s(z) \ti W^u(f^k(x)), \forall z \in \Lambda_{\epsilon} \cap B_{\nu}.$ Now, by the absolute continuity of stable lamination on Pesin's blocks, $m_{f^k(x)}^u (\mathcal{H}^s ( \Lambda_{\epsilon} \cap B_{\nu})) > 0.$ On the one hand, we know that the basin of $\nu$ is $s-$saturated so $m^u_{f^k(x)} (B ({\nu})) > 0$.  On the other hand $m^u_{f^k(x)} (B ({\mu})) =1$ and this implies that $B_{\mu} \cap B_{\nu} \neq \emptyset$ which implies $\mu=\nu$.

\end{proof}

\subsection{Uniqueness for transitive diffeomorphisms of surfaces}
Let $\mu$ and $\nu$ be two ergodic SRB measures. As mentioned before $\mu$ and $\nu$ are hyperbolic measures. By theorem \ref{periodic} we conclude that there exist hyperbolic periodic points $P_{\mu}, P_{\nu}$ such that $\nu(\Lambda(P_{\nu})) = \mu(\Lambda(P_{\mu})) = 1.$ The following proposition is the main ingredient of the proof of Theorem \ref{srb2}.

\begin{proposition} \label{homoclinic}
 $P_{\mu}$ and $P_{\nu}$ are homoclinically related and
 $ \Lambda(P_{\mu} ) = \Lambda(P_{\nu}).$
\end{proposition}

The above proposition together with  Theorem \ref{uniqueness}
immediately implies the conclusion of Theorem \ref{srb2}. To prove the existence of homoclinically relation between $P_{\mu}$ and $P_{\nu}$ we need that the manifold $M$ is two dimensional. Indeed firstly, using that $M$ is a surface we  prove that the invariant maniolds of $P_{\mu}$ and $P_{\nu}$ are topologically transverse.
Then using a finer analysis of laminations and  Sard's Theorem we prove transversal homoclinical intersection.
\begin{proof}

Recall that $P_{\mu}$ and $P_{\nu}$ (which we suppose that are
fixed points) comes from Katok's closing lemma and as $\mu, \nu$ are hyperbolic ergodic both $P_{\mu}$ and $P_{\nu}$ have non trivial homoclinic classes.  Consequently there
exist topological rectangles whose boundaries ($\partial^s$ and $\partial^u$) are consisted of
stable and unstable segments of $P_{\mu}, P_{\nu}.$ (See figue 1.)    Choose two
such rectangle $\mathcal{R}_{\mu},\mathcal{R}_{\nu}$ such that
$P_{\mu}, P_{\nu}$ respectively belong to the boundary of
$\mathcal{R}_{\mu}$ and $\mathcal{R}_{\nu}$ and $\mathcal{R}_{\mu}
\cap \mathcal{R}_{\nu} = \emptyset.$
\begin{figure}[h] \label{wow}
\psfrag{fn}{\small{$f^n(\mathcal{R}_{\mu} \cap \mathcal{R}_{\nu}) \neq \emptyset$}}
\psfrag{q}{\small{$P_{\nu}$}}
\psfrag{p}{\small{$P_{\mu}$}}
\psfrag{nu}{\small{$\mathcal{R}_{\nu}$}}
\psfrag{mus}{\small{$\mathcal{R}_{\mu}$}}
\includegraphics[width= 10cm]{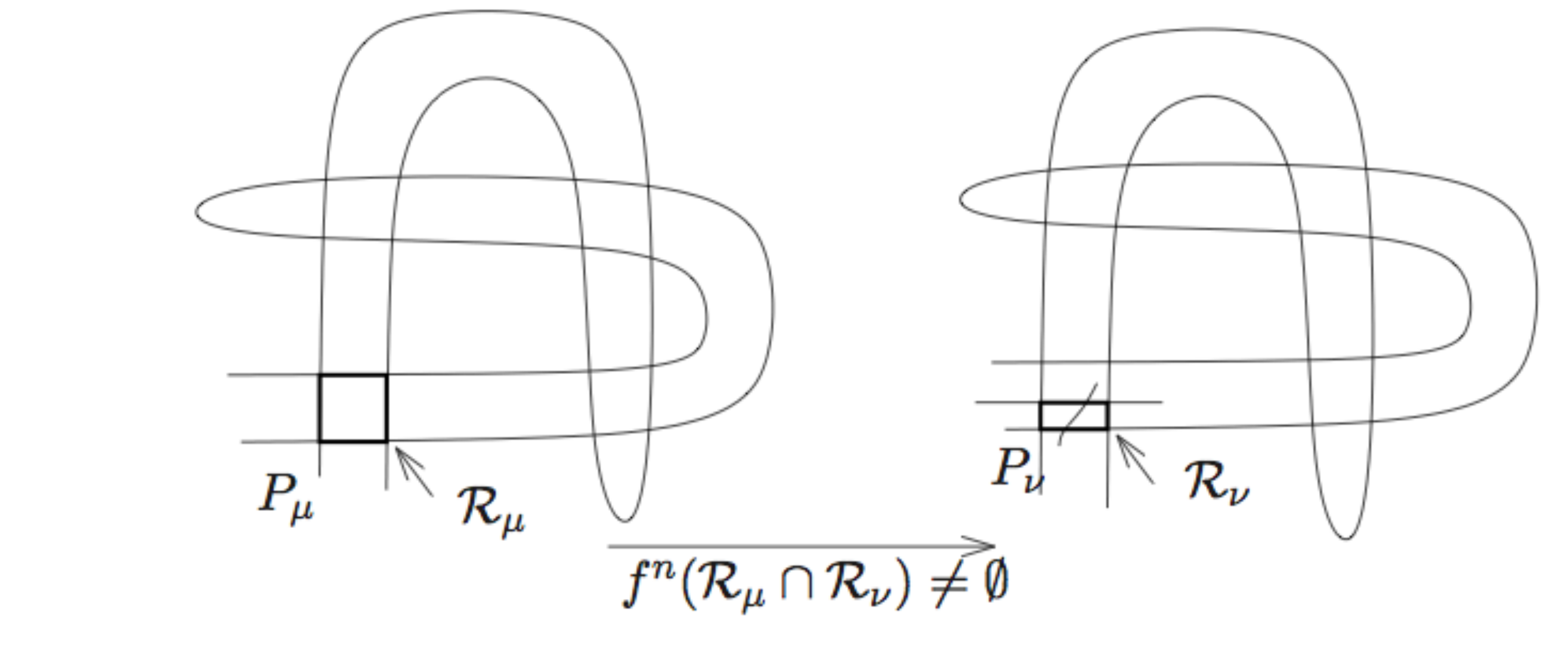}
\caption{}
\end{figure}

 By topological transitivity
of $f$ there exist $n \in \mathbb{N}$ such that
$f^n(\mathcal{R}_{\mu}) \cap \mathcal{R}_{\nu} \neq \emptyset.$
Observe that $f^n(\mathcal{R}_{\mu})$ and $\mathcal{R}_{\nu}$ are topological rectangles and
as
$P_{\nu} \notin \mathcal{R}_{\mu}$ it comes out that
$\mathcal{R}_{\nu}$ is not contained in $f^n(\mathcal{R}_{\mu}).$
So $W^u(P_{\mu}) \cap
\partial^s(\mathcal{R}_{\nu}) \neq \emptyset.$ Although $\partial^s(\mathcal{R}_{\nu})
$  is a piece of $W^s(P_{\nu})$  this intersection may be just topologically transversal (a tangency). However, we will prove that there should exist also transversal intersections between $W^u(P_{\mu})$ and $W^s(P_{\nu}).$

\begin{lemma} Taking  a rectangle $\overline{\mathcal{R}}
_{\nu}$ small enough there exist new system of coordinates such that:
 \begin{enumerate}
 \item $\overline{\mathcal{R}}_{\nu} = [0,1]^2,$
 \item   $W^u(P_{\mu}) \cap \overline{\mathcal{R}}_{\nu} $ is the graph of a $C^2$-function $\gamma : I \rightarrow [0,1], I \subset [0, 1],$
 \item There exists $K \subset [0,1]$ of positive Lebesgue measure such  $(0,x) \in [0,1]^2$ has a ``large" stable manifold crossing $[0, 1]^2.$
\end{enumerate}

\end{lemma}

\begin{proof}
We will show that it is possible to take $\overline{\mathcal{R}}_{\nu}\subset \mathcal{R}_{\nu}$ in such a way that it satisfies items 2 and 3 after a suitable change of coordinates. Observe that $W^u(P_{\mu})$ is a $C^2$-curve and we have supposed that it is tangent to local stable manifold of $P_{\nu}.$ So  to guarantee the second item of the Lemma, it is enough to take $\overline{\mathcal{R}}_{\nu}$ with small height.

To prove the last item, firstly notice that by  construction the stable (unstable) manifold of $P_{\nu}$ has transversal intersection with unstable (stable) manifold of the points inside a Pesin block of $\Lambda(P_{\nu})$.   Let us denote by  $\Lambda_{\epsilon}(P_{\nu}) $ such  Pesin hyperbolic block of $\Lambda(P_{\nu}).$ By a desintegration argument and using the fact that the conditional measures of $\mu$ along unstable manifolds are absolutely continuous with respect to Lebesgue measure  we get a local Pesin unstable manifold which intersect $\Lambda_{\epsilon}(P_{\nu})$ in a positive Lebesgue measure. Now as the stable lamination is absolutely continuous we slide such points along stable laminae to obtain a positive measure subset of $W^u(P_{\nu})$ which we denote it by $\hat{K}$.

 Now iterating $\hat{K}$ and using  $\lambda$-lemma we obtain a  positive Lebesgue measure subset $K \subset [0, 1]$ such that $\{0\} \times K \subset \partial (\mathcal{R}_{\nu})$  and $W^s(0,x)$ crosses $\overline{\mathcal{R}}_{\nu}$ for any $x \in K.$

\end{proof}

Let $\varphi(x) := h_s ( x, \gamma(x))$ where $h_s$ is the
projection by stable lamination (see Figure \ref{wow}).

\begin{figure}[h]
\psfrag{x}{$x$}\psfrag{g(x)}{$\gamma(x)$}
\psfrag{ys}{ \small{$\varphi(x) \in K$}}
\includegraphics[width= 7cm]{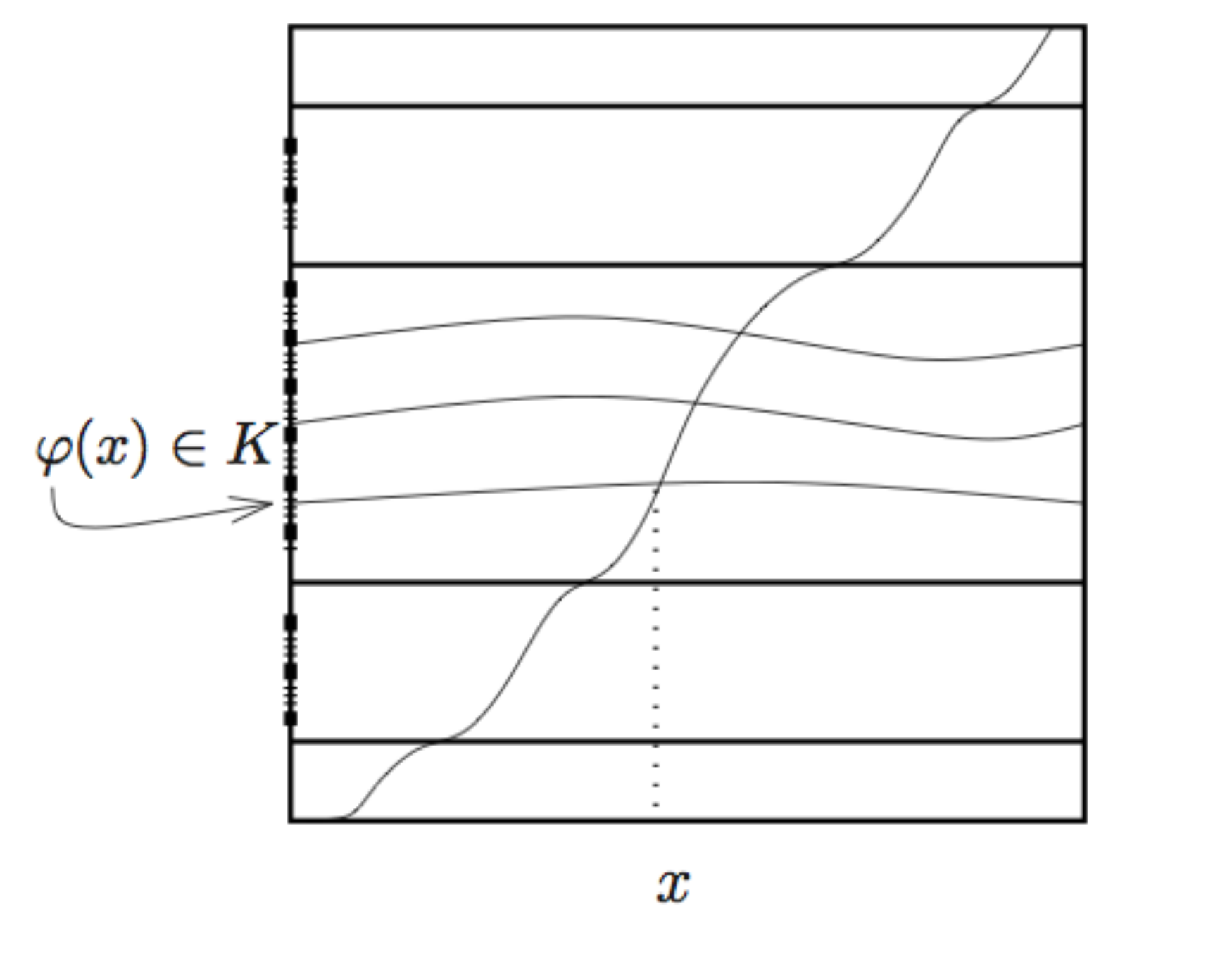}
\caption{}\label{wow}
\end{figure}

Observe that a priori $\varphi$ is
defined just on a closed subset $\{x \in I \colon (x, \gamma(x)) \in W^s(K)\}$ with positive Lebesgue measure. However we can verify the
Whitney condition to extend it in a $C^1$ fashion on the whole
inteval $I$ as we explain below.

 We recall
a standard treatment of absolute continuity of stable holonomies
in the Pesin's block following Pugh-Shub \cite{ps}. In fact we show that,
as the  stable lamination is co-dimension one its holonomy map is differentiable. More precisely we claim that the stable holonomy can be extended
to a $C^1$ function $h^s : [0, 1] \rightarrow [0, 1]$ and consequently $\varphi$ can be extended to a $C^1$-function on $I.$

Let $\mathcal{F}$ be a $C^1-$foliation which is close to the
stable lamination in $C^1$-sense. Graph transforation arguments
show that $f^{-n} (\mathcal{F})$ converges to the stable
lamination.  Let $(h_n, Jh_n)$ represent  the holonomy $h_n$  of
$f^{-n}(\mathcal{F})$ with its derivative $Jh_n$. As the domain of $h_n$ is one dimensional $Jh_n$ reprsents both the derivative and the jacobain of the holonomy $h_n$.

 By
definition the following diagram comutes
  \[
\begin{CD}
{f^{-n}(D)} @> {h_0} >> {f^{-n}(D^{'})} \\
@ V{f^{-n}} VV @ VV {f^{-n}} V \\
{D} @>> {h_n} > {D^{'}} \\
\end{CD}
\]

All holonomies $h_n$ are diffrentiable and $$Jh_n(x) =
Jf^{-n}(f^n(h_n(x))) \circ Jh_0 \circ Jf^n(x)
$$
It is standard that $(h_n, Jh_n)$ converge uniformly to $(h, Jh)$
where $h$ is the stable lamination holonomy and  for some function $Jh$.  As $h_n$ are diffrentiable with derivative $Jh_n$, by uniform convergence of $Jh_n$ to $Jh$ we conclude that $Jh$ satisfies the $C^1$-extension whitney theorem hypothesis and consequently we  can extend $\varphi$
to a $C^1$-function on the whole interval $I.$

%
Now it is easy to see that if $W^s(x, \gamma(x))$ is
tangent to the graph of $\gamma$ then $D \varphi (x) =0$. By
Sard's theorem the Lebesgue measure of critical values of
$\varphi$ is zero. From this we conclude that the graph of
$\gamma$ intersect transversally $W^s(0, x)$ for some $x \in K.$

By invariance there exist $z  \in \Lambda_{\epsilon}(P_{\nu})$ such that $W^s(z) \pitchfork W^u(P_{\mu}) \neq \emptyset.$  Now using again Katok's closing lemma we find a hyperbolic  periodic point  $\hat{P}_{\nu}$ with $W^s(\hat{P_{\nu}})$ close enough to $W^s(z)$ in $C^1$-topology so that $W^s(z) \pitchfork W^u(\hat{P_{\mu}}) \neq \emptyset$. It is clear that $P_{\nu}$ and $\hat{P}_{\nu}$ are homoclinically related. So  using $\lambda-$lemma it comes out that $W^u(P_{\mu})$ has also a transversal intersection with $W^s(P_{\nu}).$ See figure 2.

A similar argumet shows that $W^u(P_{\nu})$ has a transversal intersection with $W^s(P_{\mu})$. So we proved that $P_{\mu}$ and $P_{\nu}$ are homoclinically related.
Now using Proposition \ref{homoclinic} we obtain that $\Lambda(P_{\mu}) = \Lambda(P_{\nu}).$
\end{proof}


\begin{thebibliography}{RRRRRR}

\bibitem{abv} J.F. Alves, C. Bonatti, M. Viana, SRB measures for partially hyperbolic
systems whose central direction is mostly expanding, {\it Invent. Math.} {\bf 140} (2000),
351--398.
%
%
\bibitem{barreira-pesin} L. Barreira and Y. Pesin, Lyapunov Exponents
and Smooth Ergodic Theory, American Mathematical Society, University
Lecture Series, Vol 23.
%
\bibitem{by} M. Benedicks, L.-S. Young, Sinai-Bowen-Ruelle measure for certain H\'{e}non
maps, {\it Invent. Math.} {\bf 112} (1993), 541-576.
%
\bibitem{bv} C. Bonatti, M. Viana, SRB measures for partially hyperbolic systems whose
central direction is mostly contracting, {\it Israel J. Math.} {\bf 115} (2000), 157-194.
%
\bibitem{b} R. Bowen, {\it Equilibrium states and the ergodic theory of Anosov diffeomorphisms},
Springer Lecture Notes in Math. 470 (1975).
%

%
%
\bibitem{cy} W. Cowieson, W., L.-S. Young,
SRB measures as zero-noise limits,
{\it Ergodic Theory Dynam. Systems} {\bf 25}, (2005) 1115--1138.
%
\bibitem{d1} D. Dolgopyat, Lectures on u-Gibbs states, lecture notes, Conference on Partially
Hyperbolic Systems, Northwestern Univ. (2001).
%
\bibitem{d2} D. Dolgopyat,
On differentiability of SRB states for partially hyperbolic systems, {\it
Invent. Math.} {\bf 155}, (2004) 389--449.
%
\bibitem{furst} H. Furstenberg, Strict ergodicity and transformations of the torus, {\it Amer. J. Math.} {\bf 83}, (1961) 573--601.
    %
\bibitem{hatomoto} J. Hatomoto, Diffeomorphisms admitting SRB measures and their regularity, {\it Kodai Math. J.}  {\bf
29}, (2006) 211--226.

\bibitem{huyoung} Hu, H. Y.; Young, L.-S., Nonexistence of SBR measures for some diffeomorphisms that are ``almost Anosov, {\it Ergodic Theory Dynam. Systems} {\bf 15}, (1995), no. 1, 67--76.

%
\bibitem{kan} I. Kan, Open sets of diffeomorphisms having having two attractors, each with an everywhere dense basin, {\it Bull. Amer. Math. Soc.}  {\bf
31}, (1994) 68--74.
%
\bibitem{katok1980} A. Katok, Lyapunov exponents, entropy and
periodic orbits for diffeomorphisms, {\it IHES Publ. Math.}  {\bf
51}, (1980) 137--173.
%
\bibitem{ledrappieryoung1985} F. Ledrappier, L.-S. Young, The metric entropy of diffeomorphisms
Part I: Characterization of measures satisfying Pesin's entropy formula, {\it Ann. Math.,} {\bf 122} (1985), 509-539.
%
\bibitem{pesin1977} Ya. Pesin, Characteristic Lyapunov exponents and smooth ergodic theory, {\it Uspekhi mat. Nauk}  {\bf
32}, (1977) 55--112; English transl., {\it Russian Math. Surveys}
{\bf 32} (1977), 55--114.
%
\bibitem{psinai} Ya. Pesin,  Sinai, Ya. G.,
Gibbs measures for partially hyperbolic attractors,
{\it Ergodic Theory Dynam. Systems} {\bf 2}, (1983) 417--438.
%
\bibitem{ps} C. Pugh, M. Shub, Ergodic attractors, {\it Trans. Am. Math. Soc.} {\bf 312},  (1989) 1--54.
%
\bibitem{rhrhtu2008} F. Rodriguez Hertz, M. Rodriguez Hertz, A.
Tahzibi, R. Ures, A criterion for ergodicity of non-uniformly hyperbolic diffeomorphisms, {\it Electron. Res. Announc. Math. Sci.} {\bf 14} (2007), 74--81..
%
\bibitem{rhrhtu2009} F. Rodriguez Hertz, M. Rodriguez Hertz, A.
Tahzibi, R. Ures, New criteria for ergodicity and non-uniform hyperbolicity, preprint.
%
\bibitem{rohlin} V.A. Rohlin, On the fundamental ideas of measure theory, {\it AMS Translations} {\bf 71} (1952)
%
\bibitem{r1} D. Ruelle, A measure associated with Axiom A attractors, {\it Amer. J. Math.} {\bf 98}
(1976), 619--654.
%
\bibitem{r2} D. Ruelle, {\it Thermodynamic Formalism}, Addison Wesley, Reading, MA, 1978.
%
\bibitem{sinai} Ya. G. Sinai, Gibbs measure in ergodic theory, {\it Russian Math. Surveys} {\bf 27} (1972),
21--69.
%
\bibitem{tsujii} Tsujii, M., Physical measures for partially hyperbolic surface endomorphisms, {\it Acta Math.} {\bf 194} (2005), 37--132.
%
\bibitem{lsy_survey} Young, L.-S., What are SRB measures, and which dynamical systems have them?, {\it J. Statist. Phys.} {\bf 108} (2002), 733--754.


\end{thebibliography}
\end{document}